\documentclass[10pt,a4paper]{book}
\usepackage{amssymb,amsmath,array,amscd,makeidx,showidx,graphicx}
 \usepackage[dvips]{epsfig}
 \usepackage{amsgen, amstext,amsbsy,amsopn, amsthm, amsfonts,amssymb,amscd,amsmat
 h,euscript,enumerate,url,verbatim,calc,xypic}
 \usepackage{textcomp}
 \newcommand{\n}{\mathfrak{n} }
 \newcommand{\m}{\mathfrak{m} }
 \newcommand{\q}{\mathfrak{q} }
 \newcommand{\p}{\mathfrak{p} }

 \newcommand{\R}{\mathcal{R}}

 \newcommand{\ov}{\overline}

    \newcommand{\coker}{\operatorname{coker}}

  \newcommand{\Ass}{\operatorname{Ass}}

  \newcommand{\ann}{\operatorname{ann}}

 \newcommand{\dm}{\operatorname{dim}}
 \newcommand{\h}{\operatorname{ht}}
 
 \newcommand{\Supp}{\operatorname{Supp}}
\newcommand{\Assh}{\operatorname{Assh}}
 \newcommand{\Spec}{\operatorname{Spec}}

 \newcommand{\depth}{\operatorname{depth}}

 \newcommand{\Tor}{\operatorname{Tor}}
 \newcommand{\Hom}{\operatorname{Hom}}

  \newcommand{\lm}{\lambda}

\theoremstyle{plain}
 \newtheorem{theorem}{Theorem}[]
 
 \newtheorem{corollary}[theorem]{Corollary}
 \newtheorem{lemma}[theorem]{Lemma}
 \newtheorem{proposition}[theorem]{Proposition}

 \theoremstyle{definition}

 \newtheorem{example}[theorem]{Example}
 \theoremstyle{remark}
\begin{document}
\thispagestyle{empty}
\sloppy


\vspace{4mm}

\parindent=0mm
\begin{center}
{\Large{\bf{NEGATIVITY  OF THE CHERN NUMBER \\ [3mm]OF
PARAMETER
IDEALS}}}
\end{center}

\vspace{.3cm}
\begin{center}
{\bf{Shiro Goto$^{1}$, Mousumi Mandal$^{2}$ and Jugal Verma$^{3}$}} \\

\vspace{.3cm}
$^{1}$Department of Mathematics, School of Science and
Technology, Meiji University,\\
1-1-1 Higashi-mita, Tama-ku, Kawasaki 214-8571, Japan \\
\vspace{.3cm}
$^{{2} \& {3}}$Department of Mathematics, IIT Bombay, \\
Powai, Mumbai
400076 India\\
Email: goto@math.meiji.ac.jp, mousumi@math.iitb.ac.in,
jkv@math.iitb.ac.in
\end{center}

\markboth{Shiro Goto, Mousumi Mandal and Jugal Verma}
{{\it{Negativity of the Chern Number of Parameter Ideals}}}

\begin{center}
{\bf{1. Introduction}}
\end{center}

\markboth{Shiro Goto, Mousumi Mandal and Jugal Verma}
{{\it{Negativity of the Chern Number of Parameter Ideals}}}

Let $I$ be an ideal in a Noetherian local
 ring $R$ and let $M$ be a finite $R$-module of
dimension $d$ so that $\lm(M/IM) < \infty.$  Here $\lm(N)$
denotes the length of an $R$-module $N$. Let
$H_I(M,n)=\lm(M/I^nM)$ denote the Hilbert function of $I$
with respect to $M.$ The Hilbert function  $H_I(M, n)$ for
large $n$ is given by a polynomial $P_I(M,x)$ of degree $d.$
It is written in the form
$$P_I(M, x)=e_0(I, M){x+d-1\choose d}-e_1(I,M){x+d-2\choose
d-1}+\cdots +(-1)^de_d(I,M)$$
where $e_i(I,M)$ for $i=0,1, \ldots, d \in \mathbb Z$ are
called the Hilbert coefficients of $I$ with respect to $M.$
If $M=R$ then we write $H_I(n)=H_I(M,n), P_I(M,x)=P_I(x)$
and $e_i(I,M)=e_i(I)$ for $i=0,1, \ldots, d.$
The leading coefficient $e_0(I,M)$ is called the
{\bf multiplicity }  of $I$ with respect to $M$ and the
coefficient $e_1(I,M)$ is called the {\bf Chern number } of
$I$ with respect to $M.$
We say that  an ideal $I$ of a local ring $R$
is a {\bf parameter ideal}  for  an $R$-module $M$ of
dimension
$d$ if  $I$  is generated by $d$  elements and $\lm(M/IM) <
\infty.$ If  $I$ is a parameter ideal  of a local ring $R$
then $R$ is Cohen-Macaulay if and only if $e_0(I)=\lm(R/I)$.
Moreover if $R$ is Cohen-Macaulay then $e_1(I)=0$ for every
parameter ideal $I$  of $R.$  Vasconcelos observed
that the Chern number can also be used to characterize
Cohen-Macaulay property of $R$ for large classes  of
local rings. In the Yokohama Conference in 2008 Vasconcelos
proposed several conjectures about the Chern number of
filtrations of ideals \cite{vas}.  One of these was

\smallskip
\noindent
{\bf The Negativity Conjecture (NC):} Let $I$ be a
parameter ideal of a Noetherian local  ring $R.$
Then $e_1(I)< 0$ if and only if $R$ is not
Cohen-Macaulay.

\markboth{Shiro Goto, Mousumi Mandal and Jugal Verma}{{\it{Negativity of the Chern Number of Parameter Ideals}}}


\markboth{Shiro Goto, Mousumi Mandal and Jugal Verma}{{\it{Negativity of the Chern Number of Parameter Ideals}}}

\pagebreak

\markboth{Shiro Goto, Mousumi Mandal and Jugal Verma}{{\it{Negativity of the Chern Number of Parameter Ideals}}}


\markboth{Shiro Goto, Mousumi Mandal and Jugal Verma}{{\it{Negativity of the Chern Number of Parameter Ideals}}}

\newpage

\parindent=0mm

The objective of this survey paper is to provide a glimpse
into diverse techniques used to understand the Chern number
by presenting solutions  of NC for various classes of local
rings.

In section $1$, we present a solution of NC for
one-dimensional modules over local rings. While  doing so,
we show that the Chern number of an ideal  with respect to a
module is always non positive.

In section $2,$ using Serre's difference formula for the
Hilbert  function and the Hilbert polynomial, we   derive a
formula due to Schenzel \cite{sch}  for
all the Hilbert coefficients of standard parameter ideals
in  generalized Cohen-Macaulay local rings. The Negativity
Conjecture  for unmixed generalized Cohen-Macaulay local
rings follows from this formula.

In section $3,$ we prove  a theorem due to Vasconcelos which
asserts that if $(R,\m)$ is a Noetherian local ring of
dimension $d \geq 2$ and it embeds into a finite maximal
Cohen-Macualay $R$-module, then $R$ is Cohen-Macaulay
if and only if $e_1(I)=0$ for any parameter ideal $I$ of
$R.$ It follows as a cosequence that
the NC is true for Noetherian local domains which are
essentially of finite type over a field.

In section $4,$ the Chern number of parameter ideals in
certain quotients of regular local rings is calculated
explicitly. As a consequence a solution of NC is given
for  such rings. The proof illustrates the use of
Eagon-Northcott complex for calculation of the Chern
number.

Recall that a Noetherian local ring $(R, \m)$ is called
{\bf unmixed} if for each associated prime $\p$ of the
$\m$-adic
completion $\hat{R},$ $\dim \hat{R}/\p=\dim R.$

In section $5,$ we  present an example of a local ring
which is not unmixed in which the Chern number of some
parameter ideal
is zero. The Negativity Conjecture has recently been settled
by
Ghezzi, Goto, Hong, Phuong and Vasconcelos in \cite{gghopv}
for all unmixed local rings. We present their  solution in
section $5.$

\begin{center}
{\bf 2. Solution of NC for $1$-dimensional modules and
nonpositivity of $e_1(I,M)$}
\end{center}
In this section  we obtain  a formula for the
Chern number of a parameter ideal for  $1$-dimensional
finite modules over
a Noetherian local ring. This
generalizes a result of Goto-Nishida \cite{gni}.
We also show that  that Chern number of a parameter ideal
 in a local ring $R$ with respect to a finite $R$-module
is nonpositive.
\begin{theorem}[Goto-Nishida, \cite{gni}]
\label{p1}
 Let $(R,\m)$ be a Noetherian   ring and $M$ be a
finite $R$-module with $\dm M=1$. If $(a)$ is a parameter
ideal for
$M$ then $e_1((a),M)=-\lm(H^0_{\m}(M))$.
\end{theorem}

{\it {Proof: }} Let $N=H^0_{\m}(M)$ and $\ov M=M/N$. Notice that
$H^0_{\m}(\ov M)=0$ and $\dm \ov M=\dm M=1$, which implies
$\depth \ov M=1$. Thus $\ov M$ is Cohen-Macaulay
$R$-module.
 Consider the exact sequence
 $$0\longrightarrow N\longrightarrow M\longrightarrow \ov
M\longrightarrow 0.$$
By taking tensor product with $R/(a)^n$ we get the  exact
sequence for all $n\geq 1,$
\begin{equation}\label{eq}
 0\longrightarrow \ker \phi \longrightarrow
N/a^nN\xrightarrow{\phi} M/a^nM\longrightarrow \ov M/a^n\ov
M\longrightarrow 0.
 \end{equation}
 For $n>>0$ we have $(a)^nN \subseteq
\m^nH^0_{\m}(M)=0. $  Thus
\begin{eqnarray}
 \lm(\ov M/a^n\ov M)&=& \lm(M/a^nM+N)\\ \nonumber
                    &=& \lm(M/a^nM)-\lm
\left(\frac{a^nM+N}{a^nM} \right)\\ \nonumber
		    &=& \lm(M/a^nM)-\lm\left
(\frac{N}{a^nM\cap N}\right)\\ \nonumber
		    &=& \lm(M/a^nM)-\lm(N).\label{eq2}
\end{eqnarray}
From $(\ref{eq})$ and $(\ref{eq2}),$ we get  that for all
large $n,$
$$\lm(\ker \phi)=\lm(N)-\lm(M/a^nM)+\lm(\ov M/a^n \ov M)=0$$
which gives $\ker \phi =0$. Thus we get the following exact
sequence for large $n$,
 $$0\longrightarrow N\longrightarrow M/{a}^nM\longrightarrow
\ov M/{a}^n\ov M\longrightarrow 0.$$ Hence we have
$\lm(N)+\lm(\ov M/ a^n\ov M)=\lm(M/ a^nM)$. Since $\ov M$ is
Cohen-Macaulay, $$\lm(\ov M/ a^n\ov M)=e_0( (a^n),\ov
M)=e_0((a),\ov M)n=e_0((a),M)n.$$ Also for large $n$,
$\lm(M/ a^nM)=ne_0((a),M)-e_1((a),M)$. Therefore
$$~~~~~~~~~~~~~~~~~~~e_1((a),M)=-\lm(H^0_{\m}(M)).$$

\begin{flushright}$~\square$ \end{flushright}

\begin{proposition}
 Let $(R,\m)$ be a Noetherian  local  ring and let $M$ be a
finite $R$-module with $\dm M=1$. Let $a$ be a parameter for
$M$. Then $e_1((a),M)<0$ if and only if $M$ is not a
Cohen-Macaulay module.
\end{proposition}
{\it Proof: }
 Let $M$ be not Cohen-Macaulay. Then
$H^0_{\m}(M)\not= 0$. By Theorem \ref{p1},
$e_1((a),M)=-\lm(H^0_{\m}(M))< 0$.
The converse is well known \cite[Theorem 1.1.8]{bh}.~$\square$

\begin{theorem}[Mandal-Singh-Verma, \cite{msv}]
 Let $(R,\m)$ be a Noetherian local   ring and let $M$ be a
finite $R$-module
  with $\dm M=d$. Let $J$ be an ideal generated by  a system
of parameters for $M$. Then $$e_1(J,M)\leq 0.$$
\end{theorem}
{\it Proof: } Apply induction on $d$. The case $d=1$  is already proved.
Suppose $d=2$.
Let $J=(x,y)$  where $x,y$ is a superficial sequence for $J$
with respect to $M$.
 Consider the exact sequence
 $$0\longrightarrow M/(0:_Mx)\stackrel{x}\longrightarrow
M\longrightarrow M/xM\longrightarrow 0.$$
Applying $H^0_{\m}(.)$  we get
\begin{equation}\label{ses1}
0\longrightarrow
H^0_{\m}(M/(0:_Mx))\stackrel{x}\longrightarrow
H^0_{\m}(M)\stackrel{g}\longrightarrow
H^0_{\m}(M/xM)\longrightarrow C\longrightarrow 0
\end{equation}
where $C=\coker g$.
Consider the exact sequence
$$0\longrightarrow (0:_Mx)\longrightarrow M\longrightarrow
M/(0:_Mx)\longrightarrow 0.$$
Applying $H^0_{\m}(.)$ to the exact sequence we get
$$0\longrightarrow H^0_{\m}(0:_Mx)\longrightarrow
H^0_{\m}(M) \longrightarrow
H^0_{\m}(M/(0:_Mx))\longrightarrow 0.$$
Since $H^0_{\m}(0:_Mx)=(0:_Mx)$, we have
$$\lm(0:_Mx)=\lm(H^0_{\m}(M))-\lm(H^0_{\m}(M/(0:_Mx))).$$
Subtracting $\lm(H^0_{\m}(M/xM))$ from both sides of the
above equation we get
\begin{eqnarray*}
&&\lm(0:_Mx)-\lm(H^0_{\m}(M/xM))\\
&=&
\lm(H^0_{\m}(M))-\lm(H^0_{\m}(M/xM))
-\lm(H^0_{\m}(M/(0:_Mx))).
\end{eqnarray*}

\noindent
From the exact sequence (\ref{ses1}) we get
$$\lm(H^0_{\m}(M/(0:_Mx)))-\lm(H^0_{\m}(M))+\lm(H^0_{\m}
(M/xM))=\lm(C).$$ Therefore we have
$\lm(0:_Mx)-\lm(H^0_{\m}(M/xM))=-\lm(C)$. By a module
theoretic version of  \cite[Theorem 70] {mur}  we get
$$e_1(\ov J,\ov M)=e_1(J,M)-\lm(0:_Mx).$$
By Proposition \ref{p1}, $e_1(\ov J,\ov
M)=-\lm(H^0_{\m}(M/xM))$. Therefore
$$e_1(J,M)=\lm(0:_Mx)-\lm(H^0_{\m}(M/xM))=-\lm(C)\leq 0.$$
Let $d\geq 3$ and $a\in J$ be  superficial for $J$ with
respect to $M$. Since  $e_1(J,M)=e_1(J/(a),M/aM),$ we are
done by induction.~$\square$\\
\begin{center}
{\bf 3. Hilbert polynomial of standard system of
parameters in generalized Cohen-Macaulay local rings}
\end{center}

A local ring $(R,\m)$ of dimension $d$ is said to be generalized
Cohen-Macaulay local ring if $\lm(H^i_{\m}(R))<\infty $ for
$i=0,1,\ldots d-1$.
A system of parameters   $a_1,\ldots , a_d$ is
called a {\bf standard system of parameters} for  $R$ if for $\mathfrak q=(a_1, a_2, \ldots, a_d),$
$$\lm(R/\mathfrak q)-e_0(\mathfrak q) =\sum_{i=0}^{d-1}{d-1\choose i}\lm(H^i_{\m}(R)).$$
A local ring $R$ is called {\bf Buchsbaum} if for every 
parameter ideal 
$\mathfrak q$ of $R,$ the difference $ \lm(R/\mathfrak
q)-e_0(\mathfrak q)$ is independent of $\mathfrak q.$ In a
Buchsbaum local ring, every parameter ideal is standard.
Moreover, a Buchsbaum ring is generalized Cohen-Macaulay.

In this section we will give a new proof of a formula of Schenzel \cite{sch}  for the coefficients of the Hilbert polynomial of a  parameter ideal generated by a standard system of parameters. This is achieved via  an application of Serre's formula \cite{bh} for the difference of Hilbert
polynomial and Hilbert function of the associated graded ring $G(I)$ in terms of its local cohomology modules.
\begin{example}
Let $k$ be a field and $S=k[|x,y|]$. Let $R=S/(x)\cap
(x^3,y)$. Then $R$ is a $1$-dimensional generalized
Cohen-Macaulay but  not Buchsbaum. We prove
that $(y)$ is a standard system of parameters of $R$.  By
associativity formula we have $$e_0(y,R)=\sum_{P\in \Ass R,
\dim R=\dim R/P}e_0(y,R/P)\lm(R_P)=e_0(y,k[|y|])=1.$$
Note that $\lm(R/yR)=\lm(k[|x,y|]/(x^3,xy,y))=3$. Hence
$\lm(R/yR)-e_0(y,R)=2.$ On the other hand
$$H^0_{\m}(R)=\dfrac{[(x):\m^{\infty}]\cap
[(x^3,y):\m^{\infty}]}{(x^3,xy)}=\dfrac{(x)}{(x^3,xy)}.$$
Thus $\lm(H^0_{\m}(R))=2$. Hence $y$ is a standard system of
parameter of the generalized Cohen-Macaulay ring $R$.

\end{example}

\begin{theorem}
Let $(R,\m)$ be a generalized Cohen-Macaulay local ring and
$\mathfrak q$ be an ideal generated by a standard system of parameters of $R$. Let
$G(\mathfrak q)=\bigoplus_{n\geq 0}\mathfrak q^n/\mathfrak
q^{n+1}$ be
the associated graded ring.  Then $\lm(H^i_{M}(G(\mathfrak
q)))<\infty$ for $i=0,\ldots ,d-1$ where $M$ is the maximal
homogeneous ideal of $G(\mathfrak q)$
\end{theorem}

\begin{theorem}[Goto, \cite{g}]\label{th1}
Let $(R,\m)$ be a generalized Cohen-Macaulay local ring of
dimension $d$. Let $I$ be a standard parameter ideal and
$G(I)=\bigoplus_{n\geq 0} I^n/I^{n+1}$ be the associated
graded ring of $I$. Then
\begin{itemize}
\item[(i)] $[H^i_M(G(I))]_n=(0)~~(n\not= -i)$ and
$[H^i_M(G(I))]_{-i}=H^i_{\m}(R)$ for all $0\leq i<d$.\\
\item[(ii)] $[H^d_M(G(I))]_n=(0)$ for $n>-d$.
\end{itemize}
\end{theorem}

\begin{theorem}[Schenzel, \cite{sch}] \label{th2}
Let $(R,\m)$ be a generalized Cohen-Macaulay ring of
dimension $d$ and $I$ be a standard parameter ideal. Then for
$i=0,1, \ldots, d-1,$
$$e_{d-i}=(-1)^{d-i}\sum_{j=0}^i{i-1\choose
j-1}\lm(H^j_{\m}(R)).$$
\end{theorem}

{\it Proof: } Use induction on $i$.
By  Serre's difference formula for the graded ring
$G(I)=\bigoplus_{n\geq 0}I^n/I^{n+1}$ with the maximal
homogeneous ideal $M=\m G+G_+$:
\begin{equation}\label{eq1}
H_I^0(n)-P_I^0(n)=\sum_{i=0}^d(-1)^i \lm((H^i_{M}G(I))_n).
\end{equation}
where $H_I^0(n)=\lm(I^n/I^{n+1})$ is the Hilbert function of
$G(I)$ and $P_I^0(x)$ is the corresponding Hilbert
polynomial
written as follows
$$P_I^0(x)=e_0{x+d-1\choose d-1}-e_1{x+d-2\choose
d-2}+\cdots
+(-1)^{d-1}e_{d-1}.$$
Let us first prove the case $i=0$. Notice that by Theorem
\ref{th1} we have $H_I^0(n)=P_I^0(n)$ for $n\geq 1$. In
equation
(\ref{eq1}) substituting $n=0,1,\ldots ,n-1$ and adding both
sides we get
\begin{eqnarray*}
&&\sum_{i=0}^{n-1}\lm(I^i/I^{i+1})-\left[e_0{n+d-1\choose
d}-e_1{n+d-2\choose d-1} +\cdots
+(-1)^{d-1}ne_{d-1}\right]\\
&=&\lm(H^0_{\m}(R))
\end{eqnarray*}

which implies $(-1)^de_d=\lm(H^0_{\m}(R))$.
Next assume that the result is true for $i-1$ and we prove
it for $i$.
Substituting $n=-i$ in equation (\ref{eq1}) we get
\begin{eqnarray*}
&
&(-1)^{d-i}e_{d-i}(-1)^{i-1}+\sum_{k=1}^{i-1}(-1)^{d-i+k}e_{
d-i+k}(-1)^{i-k-1}{i-1\choose k}\\
&=&(-1)^{i-1}\lm(H^i_{\m}(R)).
\end{eqnarray*}
Substituting the
expressions for $e_{d-1},\ldots ,e_{d-i+1}$ in the above
equation to get
\begin{eqnarray*}
&&(-1)^{d-i}e_{d-i}(-1)^{i-1}+\sum_{k=1}^{i-1}(-1)^{i-k-1}{
i-1\choose k}\left[ \sum_{j=0}^{i-k}{i-k-1\choose
j-1}\lm(H^j_{\m}(R))\right]\\
&=&(-1)^{i-1}\lm(H^i_{\m}(R))
\end{eqnarray*}

From the above equation coefficient of
$\lm(H^{i-k}_{\m}(R))$ for $1\leq k\leq i$ is
\begin{eqnarray*}
&=&(-1)^{i-1}\sum_{l=1}^k(-1)^l{i-1\choose l}{i-l-1\choose
i-k-1}\\
&=&(-1)^{i-1}\sum_{l=1}^k(-1)^l{i-1\choose i-k-1}{k\choose
l}\\
&=&(-1)^i {i-1\choose i-k-1}
\end{eqnarray*}
Thus we have
\begin{eqnarray*}
& &(-1)^{d-i}e_{d-i}(-1)^{i-1}\\
&=&(-1)^{i-1}\left[\lm(H^{i}_{\m}(R))+\sum_{k=1}^{i}{i-1\choose
i-k-1}\lm(H^{i-k}_{\m}(R))\right]\\
&=&
(-1)^{i-1}\left[\lm(H^{i}_{\m}(R))+\sum_{j=0}^{i-1}{i-1\choose
j-1}\lm(H^{j}_{\m}(R))\right]
\end{eqnarray*}
Hence we get
$$(-1)^{d-i}e_{d-i}=\sum_{j=0}^{i}{i-1\choose
j-1}\lm(H^j_{\m}(R)).~\square $$

The above formula yields a characterization of depth in
terms 
of vanishing of Hilbert coefficients in unmixed generalized
Cohen-Macaulay local rings. In particular it solves the NC
for such rings.
\begin{corollary}
Let $(R,\m)$ be an unmixed  generalized Cohen-Macaulay local
ring of positive dimension $d$ and  let
$I$ be a standard parameter ideal. Fix any  $i=1, 2, \ldots,
d.$ Then
$$\depth R\geq d-i+1 \;\; \text{ if and only if}\;\; e_{i}(I)=0.$$
\end{corollary}

{\it Proof: } If  $\depth R\geq d-i+1$ then
$H^{j}_{\m}(R)=0$ for $j=0,\ldots, d-i$. Hence by Theorem
\ref{th2} we get  $e_{i}=0$. Conversely if $e_{i}=0$ then
$H^j_{\m}(R)=0$ for $j=1,\ldots ,d-i$ and since $R$ is
unmixed $H^0_{\m}(R)=0$. Thus $\depth R\geq d-i+1 $.
\begin{flushright}
$\square$
\end{flushright}


\begin{center}
{\bf 4. Solution of NC  for local  rings
of affine domains over a field. }
\end{center}
Vasconcelos took the first step for the solution of NC
\cite{vas}. He proved that  if $R$ is a non Cohen-Macaulay
Noetherian   local domain which is essentially of finite
type over a field  then  the Chern number of any  parameter
ideal in $R$ is negative. The key idea in the proof is to
embed the local domain  in finite maximal  Cohen-Macaulay
module  and use the following theorem.
\begin{theorem}
 [\bf Vasconcelos]
 \label{e1sMCM} Let
$(R,\mathfrak{m})$ be a
Noetherian   local ring of
dimension $d\geq 2$. Suppose there is an embedding of
$R$-modules:
\[ 0 \longrightarrow R \longrightarrow E \longrightarrow C
\longrightarrow 0,\] where $E$ is a finitely generated
maximal Cohen-Macaulay $R$-module and $C=E/R.$ If $R$ is
not
Cohen-Macaulay,
then $e_1(J)<0$ for  any parameter ideal $J$.
\end{theorem}
{\it Proof: } We may assume that the residue field of $R$ is infinite.
Notice that $\depth R\geq 1$. We prove the theorem  by
  induction on $d$. Let  $d=2$, let $J$ be a parameter
ideal of $R.$  If $R$ is not Cohen-Macaulay, then by depth
lemma we have $\depth C=0$.

Let $J=(x,y)$ where we may assume that  $x$ is a superficial
element for $J$ and  $x$ does not belong to any non-maximal
associated prime of $C.$
Tensoring the exact sequence above by $R/(x)$, we get the
exact sequence
\[ 0 \longrightarrow T = \Tor_1^R(R/(x), C) \longrightarrow
R/(x) \longrightarrow E/xE \longrightarrow C/xC
\longrightarrow
0,\] where $T$ is a nonzero module of finite length. Denote
by $S$ the image of $R'=R/(x)$ in $E/xE $. Note that $S$ is
a Cohen-Macaulay ring
of dimension $1$. By the Artin-Rees Lemma, for $n >>0$,
$T\cap (y^n)R'=0$, and therefore from the diagram

\[
\diagram
0 \rto  & T \cap (y^n)R' \rto\dto & (y^n)R' \rto \dto
&(y^n)S
\rto\dto  &0 \\
0 \rto          &  T \rto                 & R' \rto  & S
\rto  &0
\enddiagram
\]
using Snake Lemma,  for large $n,$
$$\lm(R'/y^nR')=\lm(T)+\lm(S/y^nS).$$ Hence comparing the
coefficients of the Hilbert polynomial from both sides we
get
\[ e_0 n - e_1 = e_0(yS)n + \lambda(T).\]
Therefore
\begin{eqnarray} \label{e1dim1}
 e_1(J) = - \lambda(T)< 0.
\end{eqnarray}
Assume now that $d\geq 3$  and let $x$ be a superficial
element for $J$, and the modules $E$ and $C$. In the exact
sequence
\begin{eqnarray} \label{Rmodx}
 0 \longrightarrow T = \Tor_1^R(R/(x), C) \longrightarrow
R'= R/(x) \longrightarrow E/xE \longrightarrow C/xC
\longrightarrow
0,
\end{eqnarray} $T$ is either zero, and we would go on with
the induction
procedure, or $T$ is a nonzero module with  finite support.

\medskip

If $T\neq 0$ in the exact sequence (\ref{Rmodx}),  we have
$e_1(JR')= e_1(JS)$. By the induction
argument, it suffices to prove that $S$ is not
Cohen-Macaulay.

We may assume that $R$ is a complete local   ring. Since
$R$ is
embedded in a maximal Cohen-Macaulay module, any associated
prime of
$R$ is an associated prime of $E$ and therefore it is
equidimensional. Consider the exact sequences
\[ 0 \longrightarrow T \longrightarrow R' \longrightarrow S
= R'/T \longrightarrow 0, \]

\[ 0 \longrightarrow R \stackrel{x}{\longrightarrow} R
\longrightarrow R' \longrightarrow 0.\]
From the first sequence, we get the exact sequence
\[ 0 \longrightarrow H_{\mathfrak{m}}^0(T)=T \longrightarrow
H_{\mathfrak{m}}^0(R')
\longrightarrow H_{\mathfrak{m}}^0(S) = 0,\]
 since $H_{\mathfrak{m}}^i(T)=0$ for $i>0$ and $S$ is
Cohen-Macaulay of dimension
 $\geq 2$; one also has
$H_{\mathfrak{m}}^1(S)=H_{\mathfrak{m}}^1(R')=0$. From the
second sequence,
since the associated primes of
$R$ have dimension $d$, $H_{\mathfrak{m}}^1(R)$ is a
finitely
generated $R$-module. Finally,
 by Nakayama Lemma $H_{\mathfrak{m}}^1(R)=0$, and therefore
 $T=H_{\mathfrak{m}}^0(R')=0$.  Hence $R/xR \cong S $. Since
$S$ is Cohen-Macaulay,  $R'$ is Cohen-Macaulay. Hence $R$ is
Cohen-Macaulay. This is a  contradiction.
\begin{flushright}
$\square$
\end{flushright}

In order to prove NC for affine local domains, we need the concept
of balanced Cohen-Macaulay module. We recall these concepts from
\cite {bh}. An $R$-module $M$ over a local ring
$R$ is called a {\bf big Cohen-Macaulay} module if there is a system
of parameters $x$ for $R$ which is $M$-regular. Hochster proved that
if $R$ contains a field then it has a big Cohen-Macaulay module. The module $M$ is called a {\bf balanced big Cohen-Macaulay module} if every system of parameters in $R$ is $M$-regular. Balanced big Cohen-Macaulay modules have many properties in common with finite modules. For example, their set of associated primes is finite. Griffith \cite[Theorem 3.1]{gr1} and \cite [Proposition 1.4]{gr2}
showed that if $R$ is a complete local domain then it has a countably
generated balanced big Cohen-Macaulay module.

Now using the above technique the following theorem
is proved by embedding $R$ into a countably generated  balanced big Cohen-Macaulay
module.

\begin{theorem} \label{e1geoMCM} Let $(R,\mathfrak{m})$ be a
Noetherian  local  domain essentially of finite
type over a
field.
 If $R$ is not Cohen-Macaulay, then $e_1(J)<0$
 for  any parameter ideal $J$.
\end{theorem}

{\it Proof: } Let $A$ be the integral closure of $R$ and $\widehat{R}$
be its
completion. Tensor the embedding $R\subset A$ to obtain the
embedding
\[ 0 \longrightarrow \widehat{R} \longrightarrow \widehat{R}
\otimes_R A=\widehat{A}. \]
From the properties of pseudo-geometric local  rings
\cite[Section 37]{nag}, $\widehat{A}$ is a reduced
semi-local
ring with a
decomposition
\[ \widehat{A} = A_1 \times \cdots \times A_r, \]
where each $A_i$ is a complete local   domain, of
dimension $\dim R$
and finite over $\widehat{R}$.
For each $A_i$ we make use of \cite[Theorem 3.1]{gr1}
and
\cite[Proposition 1.4]{gr2}
and pick a countably generated  balanced big Cohen-Macaulay
$A_i$-module and therefore $\widehat{R}$-module. Collecting
the $E_i$
we have an
embedding
\[  \widehat{R} \longrightarrow A_1 \times \cdots \times
A_r\longrightarrow  E= E_1
\oplus \cdots \oplus E_r.\]
As $E$ is a countably generated  balanced big Cohen-Macaulay
$\widehat{R}$-module, the argument above shows  if
$\widehat{R}$
is not Cohen-Macaulay then $e_1(J)=e_1(J\widehat{R})<0$.
\begin{flushright}
$\square$
\end{flushright}


Ghezzi, Hong and Vasconcelos in \cite{ghv} have  proved
the conjecture for universally catenarian domains and
domains that are homomorphic images of Cohen-Macualay rings.

\begin{theorem}
If  $(R, \m)$ is   Noetherian   domain of
dimension $d \geq  2,$ which is a homomorphic image of a
Cohen-Macaulay Noetherian  ring and if $R$
is not Cohen-Macaulay, then $e_1 (J) < 0$ for any parameter
ideal $J.$
\end{theorem}
\begin{theorem}
 Let $(R, \m)$ be a universally catenary integral domain
containing
a field. If R is not Cohen-Macaulay, then $e_1 (J) < 0$  for
any parameter ideal $J.$\\
\end{theorem}
\begin{center}
{\bf 5. Hilbert polynomial of parameter ideals in
 quotients of regular local   rings}
\end{center}
In this section we show that the Chern number is negative
for parameter ideals in certain unmixed quotients of regular
local rings and in some cases it is independent of the
choice of the parameter ideal by explicitly finding the Hilbert
polynomial of all parameter ideals.\\

L. Ghezzi, J. Hong and W. Vasconcelos \cite{ghv} calculated
the Chern number of any parameter  ideal in certain
quotients
of regular local rings of dimension four. We recall their
result first.

\begin{example} Let $(S,\m)$ be a four dimensional
regular local ring with $S/\m$ infinite. Let $P_1, P_2,
\ldots, P_r$ be a family of height two prime ideals of $S$
so that for $i \neq j, P_i+P_j$ is $\m$-primary. Put
$R=S/\cap_{i=1}^{r} P_i.$ Let
$J$ be an $\m$-primary  parameter ideal  of $R.$
Let $L= [\oplus_{i=1}^r S/P_i]/R.$ If $J \subseteq \ann L$
then $e_1(J)=-\lm(L)$ and $e_2(J)=0.$
\end{example}

\begin{lemma}
 Let $(S,\n)$ be an $r$-dimensional regular local ring. Let
$I$ be a height $h$  ideal of $S$.  Suppose
$a_1,\ldots ,a_d\in S$  such that  $(a_1+I,\ldots, a_d+I)$
is a system of parameters in $S/I$. Then $a_1,\ldots , a_d$
is a regular sequence in $S$.
\end{lemma}
\parindent=0mm
{\it Proof: } Let $J=(a_1,\ldots ,a_d)$. Then $\lm(S/I\otimes_S
S/J)<\infty $.
Hence by Serre's theorem \cite[Theorem 3, Chapter 5]{ser}
 $$\dm S/I +\dm S/J\leq \dm S=r.$$
 As $S$ is regular, it is catenary. Thus $d+r-\h J\leq
r$. Therefore $d\leq \h J\leq d$. Hence $\h J=d$ and
consequently $a_1,\ldots ,a_d$ is an $S$-regular sequence.

\begin{flushright}
$\square$
\end{flushright}

\begin{lemma}\label{lt}
 Let $(S,\n)$ be an $r$-dimensional regular Noetherian
ring. Let $I$ be a height $h$ Cohen-Macaulay ideal of $S$.
Suppose  $J=(a_1,\ldots ,a_d)$  such that  $(a_1+I,\ldots,
a_d+I)$ is a system of parameters in $S/I$.  Then for all
$j,n\geq 1$
$$\Tor_j^S(S/J^n,S/I)=0.$$
\end{lemma}
\parindent=0mm
{\it Proof: } We will apply induction on $n$. Let $n=1$. As $J$ is a
complete intersection, the Koszul complex $K(\underbar{a})$
of the sequence
$\underbar{a}=a_1, a_2,\ldots, a_d$
 $$K(\underbar{a}): 0\longrightarrow S\longrightarrow
S^d\longrightarrow S^{d\choose 2}\longrightarrow \ldots
\longrightarrow S^d\longrightarrow S\longrightarrow
S/J\longrightarrow 0$$
gives  a free resolution of $S/J$. Tensoring the above
complex  with $R:=S/I$ we get
 $$K(\underbar{a},S/I):0\longrightarrow R\longrightarrow
R^d\longrightarrow R^{d\choose 2}\longrightarrow \ldots
\longrightarrow R^d\longrightarrow R\longrightarrow
R/K\longrightarrow 0$$ which is the Koszul complex of
$JR=K$. As $K$ is  generated by an $R$-regular sequence, the
above is a free resolution of $R/K$. Hence
$\Tor_j^S(S/J,S/I)=0$ for all $j\geq 1$. Since $J$ is
generated by a regular sequence, $J^n/J^{n+1}$ is a free
$S/J$-module. Consider the exact sequence
 $$0\longrightarrow J^n/J^{n+1}\longrightarrow
S/J^{n+1}\longrightarrow S/J^n\longrightarrow 0.$$
 This gives rise to the long exact sequence
 $$\cdots \longrightarrow
\Tor_j^S(J^n/J^{n+1},S/I)\longrightarrow
\Tor_j^S(S/J^{n+1},S/I)\longrightarrow
\Tor_j^S(S/J^{n},S/I)\longrightarrow  $$
 By induction on $n$, it follows that
$\Tor_j^S(S/J^{n+1},S/I)=0$ for all $j\geq 1.$
\begin{flushright}
 $\square$
\end{flushright}

\begin{lemma}
 Let $J=(a_1,\ldots,a_d)$ be a complete intersection of
height $d$ in a regular local  ring $(R,\m)$. Let $L$
be an $R$-module of finite length. Then $\R(J)\otimes_R L$
is a finite $\R(J)$-module of dimension $d$ and
 $$\Supp(\R(J)\otimes_R L)=V(\m\R(J)).$$
\end{lemma}

\begin{lemma}\label{ko}
Let $S$ be a Noetherian  local ring, $a_1,\ldots ,a_d$ be a
regular sequence and $J=(a_1,\ldots ,a_d).$ Let $L$ be an
$S$-module of finite length.
If $J\subseteq \ann L$ then for all $n \geq 1$
$$\lm(\Tor_1(L,S/J^n))={n+d-1\choose d-1}\lm(L).$$
\end{lemma}

{\it Proof: } By \cite[Example 10]{kod}  for any $n>0$, $J^n$ is
generated by the maximal minors of the $n\times (n+d-1)$
matrix $A$ where
 \[
  A=\left( \begin{array}{ccccccccc}
            a_1 & a_2 & a_3 &\cdots & a_d & 0 & 0 & \cdots &
0\\
	    0   & a_1 & a_2 & \cdots & a_{d-1} & a_d & 0 &
\cdots & 0\\
	    0   &  0  & a_1 & \cdots & a_{d-2} & a_{d-1} &
a_d & \cdots & 0\\
	    0   &  0  &  0 & \cdots & a_1 & a_2 & a_3 &
\cdots & a_d
           \end{array}
\right)
 \]
By  Eagon-Northcott \cite[Theorem 2]{en}, the minimal free
resolution of $S/J^n$ is given by
\begin{equation}\label{eg}
0\longrightarrow S^{\beta_{d}}\longrightarrow
S^{\beta_{d-1}}\longrightarrow \cdots \longrightarrow
S^{\beta_1}\longrightarrow S\longrightarrow S/J^n
\longrightarrow 0
\end{equation}
where the Betti numbers of $S/J^n$ are given by
$$\beta_i^S(S/J^n)={n+d-1\choose d-i}{n+i-2\choose
i-1},~~~~~~~~~~~~1\leq i\leq k.$$
Taking tensor product of (\ref{eg}) with $L$ we get the
following complex
$$0\longrightarrow L^{\beta_{d}}\longrightarrow
L^{\beta_{d-1}}\longrightarrow \cdots \longrightarrow
L^{\beta_1}\longrightarrow L\longrightarrow L/J^nL
\longrightarrow 0.$$
Since $J\subseteq \ann L,$ the  maps in the above complex
are zero.
Hence
$$\lm(\Tor_1(L,S/J^n))=\beta_1\lm(L)={n+d-1\choose
d-1}\lm(L).$$
\begin{flushright}
$\square$
\end{flushright}

\begin{theorem}[Mandal-Verma, \cite{mv}]
 Let $(S,\n)$ be a regular local ring of dimension $r$
and $I_1,\ldots ,I_g$ be Cohen-Macaulay ideals of height $h$
which satisfy the condition: $I_i+I_j$ is $\n$-primary for
$i\not= j$. Let $R=S/I_1\cap \ldots \cap I_g$ and $d=\dm
R\geq 2$. Let $a_1,\ldots ,a_d\in S$ such
 that their images in $R$ form  a system of parameters. Let
$J=(a_1,\ldots ,a_d),$  $L=[\oplus_{i=1}^gS/I_i]/R$ and
$K=JR$.
Put $H_J(L,n)=\lm(J^n\otimes_R L)$ and  $P_J(L,n)$ be the
corresponding Hilbert polynomial. Then
$$P_J(L,n)=-e_1(K){n+d-2\choose
d-1}+e_2(K){n+d-3\choose d-2}- \cdots
+(-1)^de_d(K)+\lm(L).$$
If $ J\subseteq \ann L$  then
  $$P(K,n)= e_0(K) {n+d-1 \choose d}+ \lm(L){n+d-2\choose
d-1}+\cdots +n\lm(L).$$

\end{theorem}
{\it Proof: } First we show that $\lm(L) < \infty$. Consider the exact
sequence
 \begin{equation}\label{e1}
 0\longrightarrow \frac{S}{I_1\cap \cdots \cap
I_g}\longrightarrow N=\frac{S}{I_1}\oplus \ldots \oplus
\frac{S}{I_g}\longrightarrow L\longrightarrow 0.
 \end{equation}
 Let $P$ be a non-maximal prime ideal of $S$ not containing
any $I_1,\ldots ,I_g$. Then $L_P=0$. If there is an $i$ such
that $I_i\subseteq P$, then for $j\not= i$, $I_j\nsubseteq
P$. Hence
\begin{equation*}
(S/I_1\cap \ldots \cap I_g)_P=(S/I_1\oplus \ldots \oplus
S/I_g)_P=(S/I_i)_P.
\end{equation*}
Thus $L_P=0$. Hence $\Supp L=\{\n\}$. Thus $\lm(L)<\infty$.

By the depth lemma, $\depth R=1$. Thus $R$ is not
Cohen-Macaulay.
 Tensoring (\ref{e1}) with $S/J^n$  we get the exact
sequence
 $$\longrightarrow \Tor_1^S(R,S/J^n)\longrightarrow
\bigoplus_{i=1}^g\Tor_1^S(S/I_i,
S/J^n)\longrightarrow\Tor_1^S(L,S/J^n)$$
 $$\longrightarrow R/K^n\longrightarrow \bigoplus_{i=1}^g
S/(I_i,J^n)\longrightarrow L/J^nL\longrightarrow 0.$$
 By Lemma \ref{lt}, $\Tor_1^S(S/I_i,S/J^n)=0$ for all $i,n$.
For large $n$, $J^nL=0$ as $\lm(L) < \infty$. Hence for
large $n,$
 $$\lm(R/K^n)=e_0(K){n+d-1\choose d}-e_1(K){n+d-2\choose
d-1}+\cdots +(-1)^de_d(K).$$
 By (\ref{e1}) and additivity of $e_0(J,\_)$ we get
$e_0(K)=\sum_{i=1}^ge_0(J,S/I_i)$.
 Hence $$\lm(\Tor_1^S(L,S/J^n))-\lm(L)=\sum_{i=1}^d
(-1)^i e_i(K){n+d-1-i\choose d-i}.$$ From the exact sequence
 $$0\longrightarrow J^n\longrightarrow S\longrightarrow
S/J^n\longrightarrow 0$$
 we get
  $$\longrightarrow \Tor_1^S(J^n,L)\longrightarrow
\Tor_1^S(S,L)\longrightarrow\Tor_1^S(S/J^n,L)$$
  $$\longrightarrow  J^n\otimes_S L\longrightarrow S\otimes
L\longrightarrow L/J^nL\longrightarrow 0.$$
Hence for large $n$,
\begin{equation}\label{eq3}
\lm(\Tor_1^S(S/J^n,L))=\lm(J^n\otimes_SL)=
\left[\sum_{i=1}^d (-1)^i e_i(K){n+d-1-i\choose
d-i}\right]+\lm(L).
 \end{equation}
Since $\dm(\R(J)\otimes_RL)=d$ and $d\geq 2$, $e_1(K)<0$.
 If $J\subseteq\ann L$ then by Lemma \ref{ko}
$$\lm(\Tor_1(L,S/J^n))={n+d-1\choose d-1}\lm(L).$$
Substituting in (\ref{eq3}) we get
\begin{eqnarray*}
&&{n+d-1\choose
d-1}\lm(L)-\lm(L)\\&=&-e_1(K){n+d-2\choose
d-1}+e_2(K){n+d-3\choose d-2}- \cdots +(-1)^de_d(K).
\end{eqnarray*}
Using the  equation,
$$ {n+d-1 \choose d-1}= 1+
\displaystyle{\sum_{i=1}^{d-1}{n+d-i-1\choose d-i}} $$
we obtain $e_i(K)=(-1)^i\lm(L)$ for $i=1,2,\ldots, d-1$ and
$e_d(K)=0$.
\begin{flushright}
 $\square$
\end{flushright}


\begin{example}
We have computed the following example using the
software CoCoA which illustrates  the above result. Let
$S=k[x,y,z,w]$ and
$I_1=(x,y)$ and $I_2=(z,w)$. Let $R=S/I_1\cap I_2$ and
$Q=(x+z,y+w)$ be a system of parameters of $R$. Then
$e_1(Q)=-1$ and $e_2(Q)=0$. The cocoa program
which computes the example is given below.
 \begin{verbatim}
 Alias P:=$contrib/primary;
 Use S ::= Q[x,y,z,w];
 I1 :=Ideal(x,y);
 I2 :=Ideal(z,w);
 I :=Intersection(I1,I2);
 I;
 Ideal(xw,yz,yw,xz)
 Dim(S/I);
 2
 Depth(S/I);
 1
 Q :=Ideal(x+z,y+w);
 PS := P.PrimaryPoincare(I, Q); PS;
 (3-x) / (1-x)^2
 f(x)=3-x
 f'(1)=e_1(Q)=-1
 e_2(Q)=0
 \end{verbatim}
\end{example}

\begin{center}
{\bf 5. Solution of the Negativity Conjecture}
\end{center}
Recently Ghezzi, Goto, Hong, Ozeki, Phoung and Vasconcelos
\cite{gghopv} have  solved  the Negativity Conjecture in
full
generality. We sketch their solution in this section.
Recall that  for an $R$-module $M$  of finite Krull
dimension,

$$\Assh (M) =\{ \p \in \Ass (M) \mid \dim R/\p=\dim M\}.$$

\begin{lemma}[Goto-Nakamura, \cite{gn}]\label{lem4}
Let $(R,\m)$ be a complete local ring of dimension
$\geq 2$. Let $K_R$ be the canonical module of $R$ and
$S=\Hom_R(K_R,K_R)$. Let $\phi :R\longrightarrow S$ be the
map  $\phi(a)(x)=ax$ for $a\in R$ and $x\in
K_R$. If $\Ass R\subseteq \Assh R\cup \{\m \}$, then
$\ker \phi$ has  finite length and
$H^1_{\m}(R)\simeq H^0_{\m}(\coker
\phi)$. In particular, $H^1_{\m}(R)$ has finite length.

\end{lemma}

\begin{lemma} [Goto-Nakamura, \cite{gn}]\label{lem5}
 Let $R$ be a homomorphic image of a Cohen-Macaulay
local ring and assume that $\Ass R\subseteq \Assh
R\cup \{\m\}$. Then $$\mathcal F=\{\p\in \Spec R\mid \h
_R\p>1=\depth R_{\p},\p\not=\m\}$$
is a finite set.
\end{lemma}

\begin{proposition}\label{pro1}
 Let $R$ be a homomorphic image of a Cohen-Macaulay
local  ring and assume that $\Ass R\subseteq \Assh
R\cup \{\m\}$. Let $\q$ be a parameter ideal. Then there
exists a system $a_1,\ldots ,a_d$ of generators of $\q$ such
that $\Ass R/\q_i\subseteq \Assh R/\q_i\cup \{\m\}$ for all
$0\leq i\leq d$ where $q_i=(a_1, a_2, \ldots, a_i).$
\end{proposition}

\begin{theorem}[\cite{gghopv}]
 Let $(R,\m)$ be an unmixed  local   ring with $\dim
R=d>0$. Let $Q=(a_1,\ldots ,a_d)$ be a parameter ideal in
$R$ such that $e_1(Q)=0.$ Then $R$ is
Cohen-Macaulay.
\end{theorem}
{\it Proof: } Apply  induction on $d$. Let $d=1$. As $R$ is unmixed,
it is  Cohen-Macaulay.  Hence $e_1(Q)=0.$ Now
consider the  $d=2$ case. Since $R$ is unmixed we can choose
a nonzerodivisor $a_1 \in Q$ which is
a  superficial element for $Q.$
 Let $\ov R=R/a_1R$. We have $e_1(Q\ov R)=e_1(Q)=
0$. Since $\ov R$ is a $1$-dimensional,
 $e_1(Q\ov R)=-\lm(H^0_{\m}(\ov
R))=0$, which implies $e_1(Q\ov R)=0$. Hence $\bar{R}$ is
Cohen-Macaulay and so is $R.$

Assume that $d\geq 3$ and that our assertion
holds true for $d-1$. Then we can choose by Proposition
\ref{pro1}, an element $x=a_1$ superficial for the parameter
ideal $Q$ and $\Ass (R/xR) \subseteq \Assh (R/xR) \cup
\{\mathfrak m\}.$ Let $U$ be the unmixed component of $(0)$
in
$B=R/xR$. Let $(0)B=\displaystyle{\cap_{\p\in \Ass
B}\q(\p)}$. Note that
 \begin{eqnarray*}
  H^0_{\m}(B)&=&0:_B \m ^\infty\\
             &= & {\displaystyle\bigcap_{\p\in \Ass
B}}\q(\p):_B \m^\infty \\
	     &=& {\displaystyle\bigcap_{\p\in \Ass
B}}\q(\p):_B \m^\infty
\end{eqnarray*}
Since $(\q(\p):\m^\infty)=B$ if $\sqrt{\q(\p)}=\m$ and
$(\q(\p):\m^\infty)=\q(\p)$ otherwise, thus $H^0_{\m}(B)=U$.
Hence $U$ is of finite length. Notice that the
$(d-1)$-dimensional ring $B/U$ is Cohen-Macaulay by
the induction hypothesis because $B/U$ is unmixed and
$$e_1(Q(B/U))=e_1(QB)=e_1(Q)= 0.$$
Hence $H^i_{\m}(B/U)=0$ for $0\leq i \leq d-2$. Since
$H^i_{\m}(B)=H^i_{\m}(B/U)$ for all $i>0$, so
$H^i_{\m}(B)=0$ for $1\leq i\leq d-2$. Consider the
 short exact sequence
$$0\longrightarrow R \xrightarrow{x} R\longrightarrow B
\longrightarrow 0.$$
Applying the local cohomology functor on this sequence
we get the following long exact sequence of local
cohomology modules.
$$0\longrightarrow H^0_{\m}(R)\longrightarrow
H^0_{\m}(R)\longrightarrow H^0_{\m}(B)
\longrightarrow H^1_{\m}(R)\longrightarrow
H^1_{\m}(R)\longrightarrow
H^1_{\m}(B)\longrightarrow $$
Since $H^1_{\m}(B)=0$, $H^1_{\m}(R)=xH^1_{\m}(R)$. As
$H^1_{\ m}(R)$ is finitely generated, by Nakayama's Lemma
$H^1_{\m}(R)=0$. Since $R$ is unmixed $H^0_{\m}(R)=0$. From
the above long exact sequence we get $H^0_{\m}(B)=0$, which
implies that $B$ is Cohen-Macaulay. Hence $R$ is
Cohen-Macaulay as $x$ is regular.
\begin{flushright}
 $\square$
\end{flushright}

\begin{example}
 Let $S=k[|x,y,z,u,v,w|]$ and $I=(x,y)\cap (z,u,v,w)$ and
$R=S/I$.
 Then $Q=(x+z,u,y+v,w)$ is a system of parameters in $R$.
By computations
we have seen that $e_1(Q)=0$ but $R$ is not Cohen-Macaulay
as depth of $R$ is 1. So unmixed condition is necessary.\\
\\
\end{example}

\end{document}